\documentclass[12pt]{amsart}
\usepackage{amscd,amssymb,amsthm,verbatim}
\def\CD{\rm CD}
\def\DM{\rm DM}

\def\dem{\DM}
\def\Cdem{\DM(C)}
\def\Ndem{\DM(2^N)}

\def\barint{\mathop{-\mkern-19.5mu\int}}
\def\av#1{{\langle #1 \rangle}} %average

\newcommand{\diam}{\operatorname{diam}}

\newcommand{\dvol}{\operatorname{dvol}}

\newcommand{\Hess}{\operatorname{Hess}}

\newcommand{\Lip}{\operatorname{Lip}}

\newcommand{\R}{{\mathbb R}}

\newcommand{\Ric}{\operatorname{Ric}}

\newcommand{\supp}{\operatorname{supp}}

\newcommand{\DC}{{\cal D}{\cal C}}
\numberwithin{equation}{section}
\theoremstyle{plain}
\newtheorem{definition}[equation]{Definition}

\newtheorem{lemma}[equation]{Lemma}
\newtheorem{theorem}[equation]{Theorem}

\newtheorem{corollary}[equation]{Corollary}

\theoremstyle{remark}

\newtheorem{remark}[equation]{Remark}
\newtheorem{example}[equation]{Example}

\usepackage{graphicx}
\input{amssym.def}
\input{amssym}
\input{epsf}
\usepackage{a4wide}

\def\R{\mathbb R}

\def\var{\varepsilon}

\def\ov{\overline}
\def\cal{\mathcal}

\def\eps{\varepsilon}

\def\Hess{\mathop{{\rm Hess}\,}}

\def\2dr#1#2{\left. \frac{d^2}{d{#1}^2} \right |_{#2}}
\def\d2#1{\frac{d^2}{d{#1}^2}}

\def\med{\medskip}

\def\begeq{\begin{equation}}
\def\endeq{\end{equation}}
\def\begar{\begin{eqnarray}}
\def\endar{\end{eqnarray}}
\def\begar*{\begin{eqnarray*}}
\def\endar*{\end{eqnarray*}}
\def\begal{\begin{align}}
\def\endal{\end{align}}
\def\begal*{\begin{align*}}
\def\endal*{\end{align*}}

\theoremstyle{definition}

\theoremstyle{remark}

\newtheorem*{Thm*}{Theorem}
\newtheorem*{Lem*}{Lemma}
\newtheorem*{Conj*}{Conjecture}
\newtheorem*{Cor*}{Corollary}
\newtheorem*{Def*}{Definition}
\newtheorem*{Prop*}{Proposition}
\newtheorem*{Exo*}{Exercise}
\newtheorem*{Exs*}{Examples}
\newtheorem*{Ex*}{Example}
\newtheorem*{Rk*}{Remark}
\newtheorem*{Rks*}{Remarks}

\begin{document}

\title{Weak curvature conditions and functional inequalities} 

\author{John Lott}
\address{Department of Mathematics\\
University of Michigan\\
Ann Arbor, MI  48109-1109\\
USA} \email{lott@umich.edu}

\author{C\'edric Villani}
\address{UMPA\\ ENS Lyon\\
46 all\'ee d'Italie, 69364 Lyon Cedex 07\\
FRANCE} \email{cvillani@umpa.ens-lyon.fr}

\thanks{The research of the first author was 
supported by NSF grant DMS-0306242 and the Miller
Institute}
\date{September 26, 2006}

\begin{abstract} We give sufficient conditions for a measured length
space $(X, d, \nu)$ to admit local and global Poincar\'e inequalities,
along with a Sobolev inequality.
We first introduce a condition 
$\dem$ on $(X, d, \nu)$, defined in terms of transport of measures.
We show that $\dem$, together with a doubling condition on $\nu$,
implies a scale-invariant local Poincar\'e inequality.
We show that if $(X, d, \nu)$ has nonnegative $N$-Ricci curvature and
has unique minimizing geodesics between almost all pairs of points then
it satisfies $\dem$, with constant $2^N$.
The condition $\dem$ is preserved by measured Gromov-Hausdorff limits.

We then prove a Sobolev inequality for measured
length spaces with $N$-Ricci curvature bounded below by $K > 0$.
Finally we derive a sharp global Poincar\'e inequality.
\end{abstract}

\maketitle
There has been recent work on giving a good notion for a compact
measured length space $(X, d, \nu)$ 
to have a ``lower Ricci curvature bound''.
In our previous work~\cite{LV} we gave a notion of
$(X, d, \nu)$ having nonnegative $N$-Ricci curvature, where 
$N \in [1, \infty)$ is an effective dimension.  
The definition was in terms of the optimal transport of measures on $X$.
A notion was also given
of $(X, d, \nu)$ having $\infty$-Ricci curvature bounded
below by $K \in \R$; a closely related definition in this case was given
independently by Sturm \cite{Sturm2}. In a recent contribution,
Sturm has suggested a notion of $(X, d, \nu)$ having $N$-Ricci curvature 
bounded below by $K \in \R$ \cite{Sturm3}. These notions are preserved
by measured Gromov-Hausdorff limits; when specialized to
Riemannian manifolds, they coincide with classical Ricci curvature bounds.

Several results in Riemannian geometry have
been extended to these generalized settings. For example,
the Lichnerowicz inequality of Riemannian geometry implies that for a compact 
Riemannian manifold with Ricci curvature bounded below by $K > 0$,
the lowest positive eigenvalue of the Laplacian is bounded below by $K$.
In \cite{LV} we showed that this inequality extends to measured length spaces
with $\infty$-Ricci curvature bounded below by $K$, in the form of
a global Poincar\'e inequality.

When doing analysis on metric-measure spaces, a useful analytic property
is a ``local'' Poincar\'e inequality.
A metric-measure space $(X, d, \nu)$ admits a local Poincar\'e inequality 
if, roughly speaking, for each function
$f$ and each ball $B$ in $X$, the mean deviation (on $B$) of $f$
from its average value on $B$ is quantitatively controlled by the 
gradient of $f$ on a larger ball; see Definition \ref{LPIdef}
of Section \ref{LPI} for
a precise formulation. Cheeger showed that if a metric-measure space
has a doubling measure and admits a local Poincar\'e inequality then
it has remarkable extra local structure 
\cite{Cheeger}.

Cheeger and Colding showed that local Poincar\'e inequalities exist for
measured Gromov-Hausdorff limits of Riemannian manifolds with lower Ricci 
curvature bounds \cite{Cheeger-Colding (2000)}.  The method of proof
was to show that such Riemannian manifolds satisfy a certain ``segment
inequality''
\cite[Theorem 2.11]{Cheeger-Colding}
and then to show that the property of satisfying
the segment inequality is preserved under measured Gromov-Hausdorff 
limits \cite[Theorem 2.6]{Cheeger-Colding (2000)}. This then implies
the local Poincar\'e inequality.

Following on the work of Cheeger and Colding,
in the present paper we introduce a certain condition $\dem$ on
a measured length space, with $\dem$ being short for
``democratic''.
The condition $\dem$ is defined 
in terms of what we call ``dynamical democratic transference plans''.
A dynamical democratic transference plan 
is a measure on the space of all geodesics
with both endpoints in a given ball. The ``democratic'' condition is that
the
geodesics with a fixed initial point
must have their endpoints sweeping out the ball uniformly, and similarly
for the 
geodesics with a fixed endpoint. Roughly speaking, the condition $\dem$ 
says that there is a dynamical democratic transference plan 
so that a given point is not
hit too often by the geodesics.

We show that the condition $\dem$ is preserved by measured Gromov-Hausdorff
limits.
We show that $\dem$, together with a doubling condition on the measure,
implies a scale-invariant local Poincar\'e inequality.
We then show that if $(X, d, \nu)$ has nonnegative $N$-Ricci curvature
in the sense of \cite{LV}, and in addition
for almost all $(x_0, x_1) \in X \times X$
there is a unique minimal geodesic joining $x_0$ and $x_1$, then
$(X, d, \nu)$ satisfies $\dem$. Since nonnegative $N$-Ricci curvature
implies a doubling condition, it follows that $(X, d, \nu)$ admits a local
Poincar\'e inequality.
We do not know whether the condition of 
nonnegative $N$-Ricci curvature is sufficient in itself to imply a
local Poincar\'e inequality.

In the last section of the paper we prove a Sobolev inequality
for compact measured length spaces
with $N$-Ricci curvature bounded below by $K > 0$. Our definition
of $N$-Ricci curvature bounded below by $K$ is a variation on
Sturm's $\CD(K,N)$ condition \cite{Sturm3}. 
We use the Sobolev inequality to derive a global Poincar\'e
inequality.
In the case $N = \infty$, a global Poincar\'e inequality 
with constant $K$ was proven in~\cite{LV}; 
we show that when $N < \infty$, one can improve this
by a factor of $\frac{N}{N-1}$. In the Riemannian case, this is
the sharp Lichnerowicz inequality for the lowest positive eigenvalue
of the Laplacian \cite{Lichnerowicz}. 

The appendix contains a compactness theorem
for probability measures on spaces of geodesics.

After the research concerning $\dem$ was completed, we learned of 
preprints by Ohta~\cite{Ohta}, Renesse~\cite{vR}
and Sturm~\cite{Sturm3} that consider somewhat related conditions.
In~\cite{vR} a local Poincar\'e 
inequality is proved, also along the Cheeger-Colding lines, based
on a ``measure contraction property'' and
almost-everywhere unique geodesics.
The measure contraction property is also considered in
\cite{Ohta} and \cite{Sturm3}; compare with the proof of
Lemma \ref{lemJacRiem}. \\ \\
\noindent
{\bf Acknowledgement:} Many thanks are due to Xiao Zhong for
his explanations of local Poincar\'e inequalities and their possible
relation to $N$-Ricci curvature. We have also greatly benefited from
the explanations, references and enthusiastic support
provided by Herv\'e Pajot.
We thank Karl-Theodor Sturm for his comments
on this work and for sending a copy of his
preprint \cite{Sturm3}. 
We thank Luigi Ambrosio and Alessio Figalli for corrections
to earlier versions of this paper.
The first author thanks the UC-Berkeley mathematics department for
its hospitality while this research was performed. The second author also
thanks the Mathematical Forschungsinstitut Oberwolfach where part of
this work was written.

\section{Democratic couplings}

We recall some notation from \cite[Section 2]{LV}. 
Let $(X, d)$ be a compact length space; see
\cite{Burago-Burago-Ivanov (2001)} for background
material on such spaces. (Many results of the paper
extend to the locally compact case, but for simplicity we will
assume compactness.)
Let $\Gamma$ denote
the set of minimizing constant-speed
geodesics $\gamma \: : \: [0,1] \rightarrow X$,
with the time-$t$ evaluation map denoted by
$e_t \: : \: \Gamma \rightarrow X$. The endpoint map $E \: : \: \Gamma
\rightarrow X \times X$ is $E \: = \: (e_0, e_1)$.

We let $P(X)$ denote the Borel probability measures on $X$.
A {\em transference plan} $\pi \in P(X \times X)$ between $\mu_0, \mu_1
\in P(X)$ is a probability measure whose marginals are
$\mu_0$ and $\mu_1$.
The 2-Wasserstein space $P_2(X)$ is $P(X)$ equipped with 
the metric of optimal transport,
$W_2(\mu_0,\mu_1) = [\inf \int_{X \times X} 
d(x_0,x_1)^2\,d\pi(x_0, x_1)]^{1/2}$.
Here
the infimum is over transference plans between 
$\mu_0$ and $\mu_1$.
A transference plan is said to
be {\em optimal} if it achieves the infimum in the above variational
problem. When such a $\pi$ is given, we can disintegrate it
with respect to its first marginal $\mu_0$ or its second marginal
$\mu_1$. We write this in a slightly informal way:
\begeq\label{disintegrpi}
d\pi(x_0,x_1) = d\pi(x_1|x_0)\,d\mu_0(x_0)
= d\pi(x_0|x_1)\,d\mu_1(x_1).
\endeq

A {\em dynamical transference plan} consists of a transference plan
$\pi$ and a Borel measure $\Pi$ on $\Gamma$ such that
$E_* \Pi \: = \: \pi$; it is said to be optimal if $\pi$ itself is.
If $\Pi$ is a dynamical transference plan then for 
$t \in [0,1]$, we put $\mu_t \: = \: (e_t)_* \Pi$. Then
$\Pi$ is optimal if and only if
$\{\mu_t \}_{t \in [0,1]}$ is a Wasserstein geodesic
\cite[Lemma 2.4]{LV}.
Any Wasserstein geodesic arises from some optimal dynamical
transference plan in this way
\cite[Proposition 2.10]{LV}.

\begin{definition}[democratic coupling]
Given $\mu \in P(X)$,
the democratic transference plan between $\mu$ and itself is
the tensor product $\mu\otimes\mu \in P(X \times X)$. A probability measure
$\Pi\in P(\Gamma)$ is said to be a {\em dynamical democratic transference plan}
between $\mu$ and itself
if $E_*\Pi = \mu\otimes \mu$.
\end{definition}

\begin{example} \label{example}
Let $(X,d)$ be equipped with a reference measure $\nu \in 
P(X)$.  Suppose that one has 
almost-everywhere uniqueness of geodesics in the following sense :
\begeq\label{uniquegeod} 
\begin{cases} \text{For $\nu\otimes\nu$-almost all $(x_0,x_1)\in X\times X$,}
\\ 
\mbox{there is a unique geodesic } 
\gamma=\gamma_{x_0,x_1}\in \Gamma \mbox{ with } \gamma(0)=x_0 \mbox{ and }
\gamma(1)=x_1.\end{cases}
\endeq
Define $S \: : \: X \times X \rightarrow \Gamma$ measurably by
$S(x_0,x_1) \: = \: \gamma_{x_0,x_1}$. 
If $\mu$ is absolutely continuous 
with respect to $\nu$
then there is a unique dynamical democratic transference
plan between $\mu$ and itself given by
\begin{equation}
\Pi \: = \: S_*(\mu \otimes \mu).
\end{equation}
\end{example}

\begin{definition} \label{Ric'N}
A compact measured length space $(X,d,\nu)$ is a compact
length space $(X,d)$ equipped with a Borel probability
measure $\nu \in P(X)$.
Given $C > 0$, the triple $(X,d,\nu)$ is said to satisfy 
$\Cdem$ if for each ball $B=B_r(x)\subset X$ with $\nu[B] > 0$, there is
a dynamical democratic transference plan $\Pi$ from
$\mu \: = \: \frac{1_B}{\nu[B]} \: \nu$ 
to itself with the property that
if we put $\mu_t \: = \: (e_t)_*\Pi$ then 
\begin{equation} \label{TN(0)}
\int_0^1 \mu_t \,dt \leq \frac{C}{\nu[B]} \: \nu. \end{equation}
\end{definition}

We recall that a sequence $\{(X_i, d_i) \}_{i=1}^\infty$ of compact
metric spaces converges to a compact metric space $(X, d)$ in the
Gromov-Hausdorff topology if there is a sequence of
Borel maps $f_i \: : \: X_i \rightarrow X$ and a sequence of
positive numbers $\epsilon_i \rightarrow 0$ so that \\
1. For all $x_i, x_i^\prime \in X_i$,
$|d_X(f_i(x_i),f_i(x_i^\prime)) - d_{X_i}(x_i, x_i^\prime)|
\le \epsilon_i$. \\
2. For all $x \in X$ and all $i$, there is some $x_i \in X_i$
such that $d_X(f_i(x_i), x) \le \epsilon_i$.

The maps $f_i$ are called $\epsilon_i$-approximations.
If each $(X_i, d_i)$ is a length space then so is $(X,d)$.
A sequence $\{(X_i, d_i, \nu_i) \}_{i=1}^\infty$ of compact
measured length spaces converges to $(X, d, \nu)$ in the
measured Gromov-Hausdorff topology if in addition, one can
choose the $f_i$'s so that
$\lim_{i \rightarrow \infty} (f_i)_* \nu_i = \nu$ in the
weak-$*$ topology on $P(X)$.

\begin{theorem} \label{DMstab}
Suppose that  
$\{(X_i, d_i, \nu_i) \}_{i=1}^\infty$ is a sequence of compact
measured length spaces that converges to $(X, d, \nu)$ in the
measured Gromov-Hausdorff topology. Suppose that each ball $B$ in $X$
has $\nu[B] \: = \: \nu[\overline{B}]$.
If each $(X_i, d_i, \nu_i)$
satisfies $\Cdem$ then so does $(X, d, \nu)$.
\end{theorem}

\begin{proof}
Let $f_i:X_i\to X$ be a sequence of $\var_i$-approximations,
with $\var_i \rightarrow 0$, that realizes the Gromov-Hausdorff
convergence.
Let $B=B_r(x)$ be a ball in $X$ with $\nu[B]>0$.
For each $i$, choose a point $x_i \in X_i$ so that
$d_X(f_i(x_i), x) \: \le \: \eps_i$
and put $B_i=B_{r}(x_i)$. By elementary estimates,
\begin{equation}
 f_i^{-1}(B_{r-2\var_i}(x)) \subset B_i
        \subset f_i^{-1}(f_i(B_i)) \subset f_i^{-1}(B_{r+2\var_i}(x)).
\end{equation}
Combining this with the convergence of $(f_i)_*\nu_i$ to $\nu$, 
and the fact that $\nu[B]=\nu[\ov{B}]$, it is easy to
deduce that $\nu_i[B_i]\to \nu[B]$ (and in particular
$\nu_i[B_i]>0$ for $i$ large enough). A similar ``squeezing'' argument 
shows that $\int_{B_i}\varphi\circ f_i\,d\nu_i$ converges to
$\int_B \varphi\,d\nu$ for all nonnegative continuous functions $\varphi$.
As a consequence, if we put $\mu \: = \: \frac{1_B}{\nu[B]} \: \nu$ and
(for $i$ large enough) $\mu_i \: = \: \frac{1_{B_i}}{\nu_i[B_i]} \: \nu_i$
then $\lim_{i \rightarrow \infty} (f_i)_* \mu_i \: = \: \mu$ in 
the weak-$*$ topology.

For each $i$, we can introduce a dynamical democratic
transference plan $\Pi_i$ as in Definition \ref{Ric'N}, relative
to the ball $B_i$.
We write $\mu_{i,t} \: = \: (e_t)_* \Pi_i$.
From Theorem \ref{thm1}, there is a dynamical
transference plan $\Pi \in P(\Gamma(X))$, with associated
transference plan $\pi = E_* \Gamma$ and measures $\mu_t = (e_t)_* \Pi$,
so that 
$\lim_{i \rightarrow \infty} (f_i, f_i)_* \pi_i = \pi$ and
$\lim_{i \rightarrow \infty} (f_i)_* \mu_{i,t} = \mu_t$.

For any $F_1, F_2 \in C(X)$, we have
\begin{align}
\int_{X \times X} F_1(x) \: F_2(y) \: d\pi(x,y) \: & = \:
\lim_{i \rightarrow \infty} \int_{X \times X} F_1(x) \: F_2(y) \:
d((f_i,f_i)_* \pi_i)(x,y) \\
& = \:
\lim_{i \rightarrow \infty} \int_{X_i \times X_i} 
(f_i^*F_1)(x_i) \: (f_i^* F_2)(y_i) \:
d\pi_i(x_i,y_i) \notag \\
& = \:
\lim_{i \rightarrow \infty} \int_{X_i} (f_i^*F_1) \: d\mu_i 
\int_{X_i} (f_i^* F_2) \: d\mu_i \notag \\
& = \:
\lim_{i \rightarrow \infty} \int_{X} F_1 \: d(f_i)_*\mu_i 
\int_{X} F_2 \: d(f_i)_*\mu_i \notag \\
& = \:
\int_{X} F_1 \: d\mu 
\int_{X} F_2 \: d\mu. \notag
\end{align}
Thus $E_*\Pi=\pi=\mu\otimes\mu$, so $\Pi$ is still a dynamical
democratic transference plan.

It remains to check~\eqref{TN(0)}.
Let $\varphi$ be a nonnegative continuous function on $X$.
For large $i$,
we can write
\begin{equation} \int_0^1 \int_{X_i} (f_i)^*\varphi \, d\mu_{i,t} \:
\leq \: \frac{C}{\nu_i[B_i]} \: 
\int_{X_i} (f_i)^*\varphi\,d\nu_i.\end{equation}
In other words,
\begin{equation} \int_0^1 \int_{X} \varphi \, d(f_i)_*\mu_{i,t} \:
\leq \: \frac{C}{\nu_i[B_i]} \: \int_{X} \varphi\,d(f_i)_*\nu_i.\end{equation}
Passing to the limit as $i\to\infty$ gives
\begeq\label{tocheckunif} 
\int_0^1 \int_X \varphi \, d\mu_t \: \leq \:
\frac{C}{\nu[B]} \: \int_X \varphi\,d\nu.
\endeq
Since $\varphi$ is arbitrary, this proves~\eqref{TN(0)}.
\end{proof}

\section{From $\dem$ to a scale-invariant local Poincar\'e 
inequality} \label{LPI}

We first recall some notation and definitions about metric-measure spaces
$(X, d, \nu)$. If $B = B_r (x)$ is a ball in $X$ then we write
$\lambda B$ for $B_{\lambda r}(x)$. 
The measure $\nu$ is said to be {\em doubling} if
there is some $D > 0$ so that for all balls $B$,
$\nu[2B] \: \le \: D \: \nu[B]$. The constant $D$ is called the
{\em doubling constant}. An {\em upper gradient} for a function
$u \in C(X)$ is a Borel function $g \: : \: X \rightarrow
[0, \infty]$ such that for each curve
$\gamma \: : \: [0,1] \rightarrow X$ with finite length 
$L(\gamma)$ 
and constant speed,
\begin{equation}
\bigl| u(\gamma(1)) \: - \: u(\gamma(0)) \bigr| \: \le \:
L(\gamma) \: \int_0^1 \: g(\gamma(t)) \: dt.
\end{equation}
If $u$ is Lipschitz 
then an example of an upper gradient is
obtained by defining
\begin{equation} \label{ug} 
g(x) \: = \: 
\begin{cases} \displaystyle \limsup_{y\to x} \frac{|u(y)-u(x)|}{d(x,y)} \qquad
\text{if $x$ is not isolated} \\ \\
g(x)=0 \qquad\qquad\qquad\qquad \text{if $x$ is isolated}.
\end{cases}
\end{equation}

There are many forms of local Poincar\'e inequalities. 
The strongest one, in a certain sense, is as follows :
\begin{definition} \label{LPIdef}
A metric-measure space $(X, d, \nu)$ admits a local Poincar\'e
inequality if
there are constants $\lambda\geq 1$ and 
$P<\infty$ such that for all $u\in C(X)$ and $B=B_r(x)$ with
$\nu[B] > 0$,
each upper gradient $g$ of $u$ satisfies
\begeq\label{PI2}
\barint_B |u-\av{u}_B|\,d\nu \leq P r\, \barint_{\lambda B}
g \,d\nu.
\endeq
\end{definition}

Here the barred integral is the average (with respect to $\nu$) and
$\av{u}_B$ is the average of $u$ over the ball $B$. In the case of
a length space, the local
Poincar\'e inequality as formulated in Definition \ref{LPIdef}
actually implies stronger inequalities, for which we refer to
\cite[Chapters 4 and 9]{Heinonen}.  It is known that the property
of admitting a local Poincar\'e inequality is preserved under
measured Gromov-Hausdorff limits \cite{Keith,Koskela}.
(This is an extension of  the earlier result \cite[Theorem 9.6]{Cheeger}; 
Cheeger informs us that in unpublished work he also proved the extension.)

The use of a condition like $\dem$ to prove a local Poincar\'e
inequality is implicit in the work of Cheeger and Colding
\cite[Proof of Theorem 2.11]{Cheeger-Colding}. The next theorem
makes the link explicit.

\begin{theorem} \label{CdemPoinc}
If the compact measured length space
$(X,d,\nu)$ satisfies $\Cdem$ and $\nu$ is doubling, with doubling
constant $D$,
then $(X, d, \nu)$ admits a local Poincar\'e inequality
(\ref{PI2}) with $\lambda=2$ and
$P \: = \: 2CD$.
\end{theorem}

\begin{proof}
Let $x_0$ be a given point in $X$. Given $r > 0$,
write $B=B_r(x_0)$. Note that from the doubling condition,
$\nu[B] > 0$.
Put $\mu \: = \: \frac{1_B}{\nu[B]} \: \nu$.

For $y_0 \in X$, we have
\begin{equation}
u(y_0) \: - \: \av{u}_{B} \: = \: \int_X
\left( u(y_0) - u(y_1) \right) \: d\mu(y_1).
\end{equation}
Then
\begin{equation} \label{start} \barint_{B} |u-\av{u}_{B}|\,d\nu 
\: = \: \int_X \left|u(y_0)- \av{u}_{B}\right|\,d\mu(y_0)
\: \leq \: \int_{X\times X} |u(y_0)-u(y_1)|\,d\mu(y_0)\,
d\mu(y_1).
\end{equation}
Next, we estimate $|u(y_0)-u(y_1)|$ in terms of a geodesic
path $\gamma$ joining $y_0$ to $y_1$, where $y_0, y_1 \in B$.
The length of such a geodesic path is clearly less than $2r$.
Then, from the definition of an upper gradient,
\begin{equation} \label{u-u}
|u(y_0) - u(y_1)| \:  \leq \:
2r \int_0^1 g(\gamma(t))\,dt.
\end{equation}

Now let $\Pi$ be a dynamical democratic transference plan between $\mu$ and 
itself satisfying~\eqref{TN(0)}. Integrating~\eqref{u-u} against 
$\Pi$ gives, with $\mu_t \: = \: (e_t)_*\Pi$,
\begin{align}
\int_{X\times X} |u(y_0) - u(y_1)|\,d\mu(y_0)\,d\mu(y_1)
& \leq \int_\Gamma \left( 2r \int_0^1 g(\gamma(t))\,dt\right)
\,d\Pi(\gamma) \\
& = 2r \int_0^1 \left( \int_\Gamma g(\gamma(t))\,d\Pi(\gamma)\right)\,
dt \notag \\
& = 2r \int_0^1 \left( \int_\Gamma (g \circ e_t)\,d\Pi \right)\,
dt \notag \\
& = 2r \int_0^1 \left( \int_X g \, d (e_t)_*\Pi \right)\,
dt \notag \\
& = 2r \int_0^1 \int_X g \,d\mu_t \, dt. \notag
\end{align}
Combining this with~\eqref{start}, we conclude that
\begeq\label{u-uPi}
\barint_{B} |u-\av{u}_{B}|\,d\nu
\leq 2r \int_0^1 \int_X g \,d\mu_t\,dt.
\endeq
However, a geodesic joining two points in $B$ cannot leave
$2B$, so~\eqref{u-uPi} and $\Cdem$ together imply that
\begin{equation} \barint_{B} |u-\av{u}_{B}|\,d\nu
\leq \frac{2Cr}{\nu[B]} \int_{2B} g \,d\nu.\end{equation}
By the doubling property, 
$\frac{1}{\nu[B]} \: \leq \: \frac{D}{\nu[2B]}$.
The conclusion is that
\begeq
\barint_{B} |u-\av{u}_{B}|\,d\nu
\: \leq \: 2CDr  
\barint_{2B} g\,d\nu.
\endeq
This proves the theorem.
\end{proof}

\section{Nonnegative $N$-Ricci curvature and $\dem$}

In this section we show that a measured length 
space with nonnegative $N$-Ricci curvature
satisfies the condition 
$\Ndem$ as soon as geodesics are almost-everywhere unique.

We use the notion of nonnegative $N$-Ricci curvature from
\cite[Definition 5.12]{LV}. This is the same as the case $K=0$ of
Section \ref{Definitions}. We will be concerned 
here with the case $N < \infty$.

\begin{theorem} \label{thmddtp}
Suppose that a compact measured length space $(X,d,\nu)$ 
has nonnegative $N$-Ricci curvature, and that minimizing geodesics 
in $X$ are almost-everywhere unique in the sense of~\eqref{uniquegeod}.
Then $(X, d, \nu)$ satisfies $\Ndem$.
\end{theorem}

\begin{remark}
If $(X, d, \nu)$ has nonnegative $N$-Ricci curvature and
$(X, d)$ is nonbranching then minimizing geodesics in $X$ 
are almost-everywhere unique \cite{Sturm2}.
\end{remark}

Before proving Theorem~\ref{thmddtp}, we state a corollary:
\begin{corollary} \label{corddtp}
If a compact measured length space 
$(X,d,\nu)$ has nonnegative $N$-Ricci curvature
and almost-everywhere unique geodesics then it satisfies
the local Poincar\'e inequality of Definition~\ref{LPIdef} with
$\lambda\: = \: 2$ and $P\: = \: 2^{2N+1}$. More generally, if
$\{(X_i, d_i,\nu_i)\}_{i=1}^\infty$ is a sequence of compact
measured length spaces with nonnegative
$N$-Ricci curvature and almost-everywhere unique geodesics, and 
it converges in the measured Gromov-Hausdorff topology to
$(X,d,\nu)$, then $(X,d,\nu)$ satisfies the local 
Poincar\'e inequality of Definition~\ref{LPIdef} with 
$\lambda \: = \: 2$ and $P \: = \: 2^{2N+1}$.
\end{corollary}

\begin{proof}[Proof of Corollary~\ref{corddtp}]
First, \cite[Theorem 5.19]{LV} implies that
$(X,d,\nu)$ has nonnegative $N$-Ricci curvature. Then
the generalized Bishop-Gromov inequality of \cite[Theorem 5.31]{LV}
implies that $\nu[B]=\nu[\ov{B}]$ for each ball $B$
whose center belongs to the support of $\nu$. It also implies that
$\nu$ is doubling with constant $D=2^N$. Then the conclusion
follows from Theorems~\ref{DMstab}, \ref{CdemPoinc} and~\ref{thmddtp}.
\end{proof}

As preparation for the proof of Theorem \ref{thmddtp}, 
we first prove a lemma concerning optimal transport to delta functions.

\begin{lemma}\label{lemJacRiem}
Under the hypotheses of Theorem \ref{thmddtp},
let $B$ be an open
ball in $X$. Then for almost all $x_0 \in \supp(\nu)$,
the (unique) Wasserstein geodesic $\{ \mu_t \}_{t \in [0,1]}$
joining $\mu_0 \: = \: \delta_{x_0}$ to 
$\mu_1 \: = \: \frac{1_B}{\nu[B]} \ \: \nu$
can be written as $\mu_t \: = \: \rho_t\,\nu$ with
\begin{equation} \rho_t(x) \leq \frac1{t^N\nu[B]}.\end{equation}
\end{lemma}

\begin{proof}[Proof of Lemma~\ref{lemJacRiem}]
Let $\Pi$ be the (unique) optimal dynamical transference plan
giving rise to $\{ \mu_t \}_{t \in [0,1]}$.
Let $Y_0$ be the set of $x_0\in \supp(\nu)$ such that for
$\nu$-almost every $x\in X$ there is a unique geodesic joining
$x_0$ to $x$. By assumption, $Y_0$ has full $\nu$-measure.
Consider $x_0\in Y_0$. 
Given $t\in (0,1)$, $y\in X$ and $r > 0$, let $Z$ be the set of
endpoints $\gamma(1)$ of geodesics $\gamma$ with $\gamma(0) \: = \:
x_0$, $\gamma(t) \in B_r(y)$ and $\gamma(1) \in B$.
Then
\begin{align} \label{mutmuZ}
\mu_t [B_r(y)] \: & = \: ((e_t)_* \Pi)[B_r(y)] \: = \:
\Pi[e_t^{-1}(B_r(y))] \: = \: \Pi[e_1^{-1}(Z)] \: = \:
((e_1)_* \Pi)[Z] \\
& = \: \mu_1[Z] \: = \: \frac{\nu[Z]}{\nu[B]}.   \notag
\end{align}

If $\nu[Z] \: = \: 0$ then $\mu_t [B_r(y)] \: = \: 0$.
Otherwise,
put $\mu_1^\prime \: =\: \frac{1_Z}{\nu[Z]} \: \nu$ and let
$\{\mu^\prime_t \}_{t \in [0,1]}$ be the (unique) Wasserstein
geodesic joining $\mu^\prime_0 \: = \: \delta_{x_0}$ to $\mu_1^\prime$.
By the construction of $Z$, $\mu^\prime_t [B_r(y)] \: = \: 1$.

Put
\begin{equation} \phi(s) \: = \: \int_X (\rho'_s)^{1-1/N}\,d\nu
\end{equation}
where $\rho'_s$ is the density in the absolutely continuous part
of the Lebesgue decomposition of $\mu'_s$ with respect to $\nu$.
As $(X, d, \nu)$ has nonnegative $N$-Ricci curvature, $-\phi$
satisfies a convexity inequality on $[0,1]$.
(We use here the uniqueness of the Wasserstein geodesic 
$\{\mu^\prime_t \}_{t \in [0,1]}$. From the definition of
nonnegative $N$-Ricci curvature, {\em a priori} one only
has convexity along {\em some} Wasserstein geodesic
from $\mu^\prime_0$ to $\mu^\prime_1$.)

As $\phi(0) \: = \: 0$ and $\phi(1) \: = \: \nu[Z]^{\frac{1}{N}}$,
we obtain 
\begin{equation} \label{combine}
\phi(t) \: \geq \: t \: \nu[Z]^{\frac{1}{N}}.
\end{equation}

On the other hand, by Jensen's inequality
\begin{align} \phi(t) & = \nu[B_r(y)] \left(
\frac1{\nu[B_r(y)]} \int_{B_r(y)} 
(\rho'_t)^{1-\frac1N} \,d\nu \right)
\leq \nu[B_r(y)] \left( \frac1{\nu[B_r(y)]} 
\int_{B_r(y)}\rho'_t\,d\nu\right)^{1-\frac1N} \\
& \leq \nu[B_r(y)]^{\frac1N} \: \mu_t^\prime[B_r(y)]^{1-\frac{1}{N}} \: =
\: \nu[B_r(y)]^{\frac1N}. \notag
\end{align}
This, combined with (\ref{combine}), gives
\begin{equation} t \: \nu[Z]^{\frac1N} \leq \nu[B_r(y)]^\frac1N,
\end{equation}
Then by~\eqref{mutmuZ},
\begin{equation} 
\frac{\mu_t[B_r(y)]}{\nu[B_r(y)]} \: = \:
\frac{\frac{\nu[Z]}{\nu[B_r(y)]}}{\nu[B]}
\: \leq \: \frac{1}{t^N \: \nu[B]}.
\end{equation}
Since this is true for any ball centered at any $y$ 
and since balls generate the Borel $\sigma$-algebra,
we deduce that $\mu_t\leq \frac{\nu}{t^N\,\nu[B]}$; so
$\mu_t$ is absolutely continuous with respect to $\nu$
and its density is bounded above by $\frac{1}{t^N \: \nu[B]}$.
\end{proof}

\begin{proof}[Proof of Theorem~\ref{thmddtp}]
As in Definition \ref{Ric'N}, let $B$ be a ball in $X$ with
$\nu[B] > 0$ and
put $\mu \: = \: \frac{1_B}{\nu[B]} \: \nu$. 
As in Example \ref{example},
there is a unique
dynamical democratic transference plan from $\mu$ to itself.
We want to show that the condition of Definition 
\ref{Ric'N} is satisfied.

Define $\mu^{x_0}_t$ as in Lemma~\ref{lemJacRiem},
with density $\rho^{x_0}_t$. From Lemma \ref{lemJacRiem},
$\rho^{x_0}_t \: \leq \: \frac{1}{t^N \, \nu[B]}$.
The key point is that this is independent of $x_0$.

We now want to integrate with respect to $x_0$.
With $\mu_t$ as in Definition \ref{Ric'N} and $\varphi \in C(X)$, we have
\begin{align}
\int_X \varphi\, d\mu_t & = 
\int_X \varphi\, d(e_t)_* \Pi =
\int_X (\varphi \circ e_t)\,d\Pi =
\int_\Gamma \varphi(\gamma(t))\,d\Pi(\gamma) \\
&  =
\int_{X \times X} \varphi(\gamma_{x_0,x_1}(t))\,
d\mu(x_0)\,d\mu(x_1) \notag
\end{align}
and
\begin{equation} 
\int_X \varphi\, d\mu_t^{x_0} \: = \:
\int_X \varphi(\gamma_{x_0,x_1}(t))\, d\mu(x_1).\end{equation}
These equations show that
\begin{equation} \mu_t = \int_X \mu_t^{x_0}\,d\mu(x_0).\end{equation}
In particular, $\mu_t$ admits a density $\rho_t$, which satisfies the equation
\begin{equation} \rho_t(x) = \int_X \rho_t^{x_0}(x)\,d\mu(x_0).\end{equation}
It follows immediately that
\begin{equation} \rho_t(x) \leq \frac1{t^N\,\nu[B]}.\end{equation}

As geodesics are almost-everywhere unique, we can
apply the preceding arguments symmetrically with respect to 
the change $t \rightarrow 1-t$.
This gives
\begin{equation} \rho_t(x) \leq \frac1{(1-t)^N\,\nu[B]}.\end{equation}
Then
\begin{equation} \label{improved}
\rho_t(x) \: \leq \: \min \left( \frac{1}{t^N}, \frac1{(1-t)^N} \right)
\: \frac{1}{\nu[B]} \: \le \: \frac{2^N}{\nu[B]}.
\end{equation}
The theorem follows.
\end{proof}

\begin{remark} The above bounds~\eqref{improved} can be improved as follows.
Let $\mu= \rho\, \nu$ be a measure that is absolutely continuous
with respect to $\nu$, and otherwise arbitrary. Then there exists 
a probability measure $\Pi\in P(\Gamma)$, with $E_*\Pi=\mu\otimes\mu$,
such that $\mu_t=(e_t)_*\Pi$ admits a density $\rho_t$ with respect
to $\nu$, and
\begeq \label{condLp}
\|\rho_t\|_{L^p} \leq \min
                \left( \frac{1}{t^{N/p'}}, \frac{1}{(1-t)^{N/p'}}
                        \right)
 \|\rho\|_{L^p}
\endeq
for all $p \in (1, \infty)$,
where $p'=p/(p-1)$ is the conjugate exponent to $p$
and $\|\rho\|_{L^p}=(\int_X \rho^p\,d\nu)^{1/p}$.
Condition~\eqref{condLp}
is also stable by measured Gromov-Hausdorff limits.
Yet we prefer to focus on Condition~$\dem$ because
it is a priori weaker, and still implies the local Poincar\'e inequality.
\end{remark} 

\section{Definition of $N$-Ricci
curvature bounded below by $K$} \label{Definitions}

We recall some more notation.

For $N \in [1, \infty)$, 
the class $\DC_N$ is the set of continuous
convex functions $U \: : \: [0, \infty) \rightarrow \R$, 
with $U(0)=0$, such that the function
\begin{equation} \label{psi}
\psi(\lambda) \: = \: \lambda^N \: U \left( \lambda^{-N} \right)
\end{equation}
is convex on $(0, \infty)$.
For $N = \infty$,  
the class $\DC_\infty$ is the set of continuous
convex functions $U \: : \: [0, \infty) \rightarrow \R$, 
with $U(0)=0$, such that the function
\begin{equation} \label{psi2}
\psi(\lambda) \: = \: e^\lambda \: U \left( e^{-\lambda} \right)
\end{equation}
is convex on $(-\infty, \infty)$.
In both cases, such a $\psi$ is automatically nonincreasing
by the 
convexity of $U$. We write $U'(\infty) \: = \: \lim_{r \rightarrow \infty}
\frac{U(r)}{r}$.
If a reference probability measure $\nu \in P(X)$ is given, we define
a function $U_\nu \: : \: P(X) \rightarrow \R \cup \{ \infty \}$ by
\begin{equation} 
U_\nu(\mu) = \int_X U (\rho)\,d\nu \ + \ U'(\infty) \, \mu_s(X),\end{equation}
where $\mu \: =  \: \rho \, \nu + \mu_s$ is the Lebesgue decomposition of
$\mu$ with respect to $\nu$.

We now introduce some expressions
that played a prominent role in~\cite{CEMS01} and in~\cite{Sturm3}.
Given $K\in\R$ and $N \in (1, \infty]$, define
\begeq\label{beta}
\beta_t(x_0,x_1) \: = \: 
\begin{cases} 
e^{\frac16 K\: (1-t^2) \: d(x_0,x_1)^2} \qquad\qquad\qquad & \text{if 
$N = \infty$},\\
\infty \qquad\qquad\qquad & \text{if $N < \infty$, $K>0$ and $\alpha>\pi$},\\
\left(\frac{\sin(t\alpha)}{t\sin\alpha}\right)^{N-1}\qquad
& \text{if $N < \infty$, $K>0$ and $\alpha\in[0,\pi]$}, \\
1 \qquad \qquad\qquad & \text{if $N < \infty$ and $K=0$},\\
\left(\frac{\sinh(t\alpha)}{t\sinh\alpha}\right)^{N-1}\qquad
& \text{if $N < \infty$ and $K<0$},
\end{cases}\endeq
where
\begeq \label{alpha}
\alpha \: = \: \sqrt{\frac{|K|}{N-1}}\, d(x_0,x_1).
\endeq
When $N=1$, define 
\begin{equation}
\beta_t(x_0,x_1) = \begin{cases}
\infty \qquad \text{if $K>0$},\\
1 \qquad\quad \text{if $K\leq 0$}.
\end{cases}
\end{equation}
Although we may not write it explicitly, $\alpha$ and $\beta$ depend
on $K$ and $N$.
\med

\begin{definition} \label{N-Riccidef}
We say that $(X, d, \nu)$ has $N$-Ricci curvature bounded below by $K$
if the following condition is satisfied. Given $\mu_0, \mu_1 \in P(X)$
with support in $\supp(\nu)$, write
their Lebesgue decompositions with respect to $\nu$ as
$\mu_0=\rho_0\,\nu + \mu_{0,s}$ and $\mu_1=\rho_1\,\nu + \mu_{1,s}$,
respectively.
Then there is {\em some} 
optimal dynamical transference plan $\Pi$ from $\mu_0$ to $\mu_1$,
with corresponding Wasserstein geodesic
$\mu_t \: = \: (e_t)_*\Pi$, so that for all $U\in\DC_N$ and all $t\in [0,1]$,
we have
\begin{align} \label{ineqDCN}
U_\nu(\mu_t) \: \leq \: & (1-t) \int_{X \times X} \beta_{1-t}(x_0,x_1) \,
U \left( \frac{\rho_0(x_0)}{\beta_{1-t}(x_0,x_1)} \right)
\, d\pi(x_1|x_0)\,d\nu(x_0) \: +  \\
& t \int_{X \times X} \beta_t(x_0,x_1)\, U \left( \frac{\rho_1(x_1)}
{\beta_{t}(x_0,x_1)} \right) \, d\pi(x_0|x_1)\,d\nu(x_1) \: + \notag \\
& U'(\infty) \bigl( (1-t) \mu_{0,s}[X] \: + \: t \mu_{1,s}[X] \bigr). \notag
\end{align}
\end{definition}

Here if $\beta_t(x_0, x_1) = \infty$ then we interpret
$\beta_t(x_0,x_1)\, U \left( \frac{\rho_1(x_1)}
{\beta_{t}(x_0,x_1)} \right)$ as
$U^\prime(0) \: \rho_1(x_1)$, and similarly
$\beta_{1-t}(x_0,x_1) \,
U \left( \frac{\rho_0(x_0)}{\beta_{1-t}(x_0,x_1)} \right)$
as $U'(0)\:\rho_0(x_0)$. 
It is not difficult to show that if $N<\infty$ and $(X,d,\nu)$ 
has $N$-Ricci curvature bounded below by $K>0$, then
the diameter of the support of $\nu$ is bounded above by $\pi\sqrt{(N-1)/K}$.
In that case, the quantity $\alpha$ defined in~\eqref{alpha} 
will vary only in $[0,\pi]$ as $x_0$, $x_1$ vary in the support of $\nu$. 

\begin{remark} If $\mu_0$ and $\mu_1$ are absolutely continuous
with respect to $\nu$ then
the inequality can be rewritten in the more symmetric
form
\begin{align} \label{symform}
U_\nu(\mu_t) \: \leq \: & (1-t) \int_{X \times X} 
\frac{\beta_{1-t}(x_0,x_1)}{\rho_0(x_0)} \,
U \left(\frac{\rho_0(x_0)}{\beta_{1-t}(x_0,x_1)} \right)\,d\pi(x_0, x_1)
\: + \\
& t \int_{X \times X} \frac{\beta_t(x_0,x_1)}{\rho_1(x_1)}\, 
U \left( \frac{\rho_1(x_1)}
{\beta_{t}(x_0,x_1)} \right) \, d\pi(x_0, x_1). \notag
\end{align}
\end{remark}

\begin{remark}
Note that (\ref{ineqDCN}) is unchanged by the addition
of a linear function $r \rightarrow cr$ to $U$.
Of course, the validity of~\eqref{ineqDCN} depends on
the values of $K$ and $N$. The parameter
$\beta_t$ is monotonically nondecreasing in $K$ and the function
$\beta\longmapsto \beta U(\rho/\beta)$ is monotonically
nonincreasing in $\beta$
(because of the convexity of $U$). It follows that if $K \le K^\prime$ and
$(X, d, \nu)$ has $N$-Ricci curvature bounded below by $K^\prime$ then
it also has $N$-Ricci curvature bounded below by $K$, as one would
expect. One can also show that if $N \le N^\prime$ and
$(X, d, \nu)$ has $N$-Ricci curvature bounded below by $K$ then
it has $N^\prime$-Ricci curvature bounded below by $K$.
\end{remark}

We now compare Definition \ref{ineqDCN} with earlier definitions in
the literature, starting with the case
$N < \infty$.  If $N < \infty$ and $K=0$ then one recovers
the $N<\infty$, $K=0$ definition of \cite{LV}. If
$N < \infty$ and one specializes to $U(r)$ being 
\begin{equation} \label{UN}
U_N(r) \: = \: Nr 
\left(1 \: - \: r^{-1/N} \right),
\end{equation}
with corresponding entropy function
\begin{equation} \label{HN}
H_{N,\nu}(\mu) = N - N \int_X \rho^{1 - \frac{1}{N}} \: d\nu,
\end{equation}
then one recovers the $N<\infty$ definition of
Sturm~\cite{Sturm3}.
(In \cite{Sturm3} it was not required that $\pi$ and $\{\mu_t\}_{t \in [0,1]}$
be related in the sense that they both arise from an optimal
dynamical transference plan $\Pi$. We can make that requirement without
loss of consistency.)

To deal with the $N = \infty$ case, we use the following lemma.

\begin{lemma} \label{comp}
If $N = \infty$, with $\psi$ defined as in~\eqref{psi2}, put
\begin{equation}
\lambda(U) = 
\begin{cases}
- \: K \: \psi^\prime(\infty) & \text{ if $K > 0$}, \\
0 & \text{ if $K = 0$}, \\
- \: K \: \psi^\prime(- \infty) & \text{ if $K > 0$}.
\end{cases}
\end{equation}
If $\mu_0$ and $\mu_1$ are absolutely continuous
with respect to $\nu$ then
\begin{equation}
\int_{X \times X} \frac{\beta_t(x_0,x_1)}{\rho_1(x_1)}\, 
U \left( \frac{\rho_1(x_1)}
{\beta_{t}(x_0,x_1)} \right) \, d\pi(x_0, x_1) \: \le \:
\int_X U(\rho_1) \: d\nu \: - \: \frac{1}{6} \: (1-t^2) \: \lambda(U)
\: W_2(\mu_0, \mu_1)^2.
\end{equation}
\end{lemma}
\begin{proof}
Suppose first that $K > 0$. From the convexity of $\psi$, if
$\rho_1(x_1) > 0$ then
\begin{equation}
\frac{
\psi(- \ln \rho_1(x_1) + \frac16\, K (1-t^2) \: d(x_0, x_1)^2) \: - \:
\psi(- \ln \rho_1(x_1))}{
\frac16 K (1-t^2) \: d(x_0, x_1)^2
} \: \le \: \psi^\prime(\infty).
\end{equation} 
Then
\begin{equation}
\frac{\beta_t(x_0,x_1)}{\rho_1(x_1)}\, 
U \left( \frac{\rho_1(x_1)}
{\beta_{t}(x_0,x_1)} \right) \: \le \:
\frac{1}{\rho_1(x_1)}\, 
U \left({\rho_1(x_1)} \right)
\: + \: \frac16 K \psi^\prime(\infty) (1-t^2) \: d(x_0, x_1)^2.
\end{equation}
The lemma follows upon integration with respect to $d\pi(x_0, x_1)$.
The cases $K=0$ and $K<0$ are similar.
\end{proof}

Using Lemma \ref{comp}, and the analogous inequality for
$\int_{X \times X} 
\frac{\beta_{1-t}(x_0,x_1)}{\rho_0(x_0)} \,
U \left(\frac{\rho_0(x_0)}{\beta_{1-t}(x_0,x_1)} \right)\,d\pi(x_0, x_1)$,
one finds that (\ref{symform}) implies
\begin{equation} \label{lvineq}
U_\nu(\mu_t) \: \le \: t \: U_\nu(\mu_1) \: + \:
(1-t) \: U_\nu(\mu_0) \: - \: \frac12 \: \lambda(U) \: t(1-t) \: 
W_2(\mu_0, \mu_1)^2,
\end{equation}
which is exactly the inequality used in \cite{LV} to define what it means
for $(X,d,\nu)$ to have
$\infty$-Ricci curvature bounded below by $K$, in the sense of
\cite{LV}. We have only
shown that (\ref{lvineq}) holds when $\mu_0$ and $\mu_1$ are absolutely
continuous with respect to $\nu$, but 
\cite[Proposition 3.21]{LV} then implies that it holds for all
$\mu_0, \mu_1$ with support in $\supp(\nu)$.

A consequence is that any result of \cite{LV} concerning 
measured length spaces with $\infty$-Ricci curvature bounded below by $K$,
in the sense of \cite{LV}, also holds for measured length spaces with
$\infty$-Ricci curvature bounded below by $K$ in the sense of
Definition \ref{N-Riccidef}.

Finally, Sturm's notion of having $\infty$-Ricci curvature bounded below by $K$
\cite{Sturm2} is the specialization of the definition of
\cite{LV} to the case $U(r) = U_\infty(r) = r \ln (r)$.

The notion of having $N$-Ricci curvature bounded below by $K$,
in the sense of Definition \ref{N-Riccidef}, is preserved by
measured Gromov-Hausdorff limits.  We will present the proof,
which is more complicated than that of the analogous statement in
\cite{LV}, elsewhere. 

We now show that in the case of a Riemannian manifold 
equipped with a smooth
measure, a lower Ricci curvature bound in the
sense of Definition \ref{N-Riccidef} is equivalent to a tensor inequality.

Let $M$ be
a smooth compact connected $n$-dimensional manifold with Riemannian
metric $g$. We let $(M, g)$ denote the corresponding metric space. 
Given $\Psi \in C^\infty(M)$ with $\int_M e^{-\Psi} \: \dvol_M \: = \: 1$,
put $d\nu \: = \: e^{-\Psi} \: \dvol_M$.

\begin{definition}
For $N \in [1, \infty]$, let the $N$-Ricci tensor $\Ric_N$ of $(M, g, \nu)$
be defined by
\begin{equation}
\Ric_N \: = 
\begin{cases}
\Ric \: + \: \Hess(\Psi)  & \text{ if $N = \infty$}, \\
\Ric \: + \: \Hess(\Psi) \: - \: \frac{1}{N-n} \: d \Psi \otimes
d \Psi & \text{ if $n \: < \: N \: < \: \infty$}, \\
\Ric \: + \: \Hess(\Psi) \: - \: \infty \: (d \Psi \otimes
d \Psi) & \text{ if $N = n$}, \\
- \infty & \text{ if $N < n$,}
\end{cases}
\end{equation}
where by convention $\infty \cdot 0 \: = \: 0$.
\end{definition}

\begin{theorem} \label{Riemequiv}
For $N \in [1, \infty]$, the measured length space 
$(M, g, \nu)$ has $N$-Ricci curvature bounded below by $K$ if and
only if $\Ric_N \: \ge \: Kg$.
\end{theorem}
\begin{proof}
The proof is similar to that of \cite[Theorems 7.3 and 7.42]{LV} and
\cite[Theorem 1.9]{Sturm3}. We only sketch a few points of the
proof, in order to clarify the role played by the function $U$.

Suppose that $N \in (1, \infty)$ and
$\Ric_N \: \ge \: Kg$. We want to show that the condition
in Definition \ref{N-Riccidef} holds.
As in \cite[Theorem 7.3]{LV}, 
we can reduce to the case when $\mu_0$ and $\mu_1$
are absolutely continuous with respect to $\nu$. The unique
Wasserstein geodesic between them is of the form $\mu_t \: = \: (F_t)_*
\mu_0$ for certain maps $F_t \: : \: M \rightarrow M$. Put
\begin{equation}
C(y,t) 
\: = \: e^{- \: \frac{\Psi(F_t(y))}{N}} \: {\det}^{\frac{1}{N}}(dF_t)(y)
\end{equation} and
\begin{equation}
\eta_0 \: = \: \frac{d\mu_0}{\dvol_M}.
\end{equation}
Then in terms of the function $\psi$ of (\ref{psi})
there is an equation \cite[(7.19)]{LV}
\begin{equation} \label{touse}
U_{\nu}(\mu_t) \:
 = \:
\int_M  \psi \left( C(y, t) \: \eta_0^{-\frac{1}{N}}(y) \right) \: d\mu_0(y).
\end{equation}

With the notation~\eqref{alpha} in use, define
\begin{equation} \tau_{K,N}^{(t)} (d(x_0,x_1))\:= 
\begin{cases} t^{\frac1N} \left(\frac{
\sin (t\alpha)}{\sin\alpha}\right)^{1-\frac1N} & \text{ if } K>0, \\
t & \text{ if } K=0,\\
t^{\frac1N} \left(\frac{
\sinh (t\alpha)}{\sinh\alpha}\right)^{1-\frac1N} & \text{ if } K<0,
\end{cases}
\end{equation}
one can show by combining \cite[Section 7]{LV} and
\cite[Section 5]{Sturm3} that
\begin{equation}
C(y,t) \: \ge \: 
\tau_{K,N}^{(1-t)}(d(y,F_1(y)))\,
C(y, 0) \: + \:
\tau_{K,N}^{(t)} (d(y,F_1(y)))\, C(y,1).
\end{equation}
As $\psi$ is nonincreasing, we obtain
\begin{equation}
U_{\nu}(\mu_t) \: 
 \le \:
\int_M \psi\left( \frac{\tau_{K,N}^{(1-t)}(d(y,F_1(y)))\,
C(y, 0) \: + \:
\tau_{K,N}^{(t)} (d(y,F_1(y)))\, C(y,1)}{\eta_0^{\frac1N}(y)}
\right)\,d\mu_0(y).
\end{equation}
As $\psi$ is convex by assumption, we now obtain
\begin{align}
U_{\nu}(\mu_t) \: \le \:
& (1-t) 
\int_M \psi\left(\frac{\tau_{K,N}^{(1-t)}(d(y,F_1(y)))}{1-t}\,
\frac{C(y,0)}{\eta_0^{\frac1N}(y)}\right) \: d\mu_0(y) \: + \\
& t \int_M \psi\left(\frac{\tau_{K,N}^{(t)}(d(y,F_1(y)))}{t}\,
\frac{C(y,1)}{\eta_0^{\frac1N}(y)}\right) \: d\mu_0(y).\notag
\end{align}
After using the definition of $\psi$ again, along with
(\ref{touse}) in the cases $t=0$ and $t=1$,
one arrives at~\eqref{symform}.
\end{proof}

The next result is an analog of \cite[Theorem 5.52]{LV}.

\begin{theorem} \label{ac}
If $(X, d, \nu)$ has $N$-Ricci curvature bounded below by $K$ then
for any $\mu_0, \mu_1 \in P(X)$ that are absolutely continuous with
respect to $\nu$, the Wasserstein geodesic $\{\mu_t\}_{t \in [0,1]}$
of Definition \ref{N-Riccidef} has the property that $\mu_t$ is absolutely
continuous with respect to $\nu$, for all $t \in [0,1]$.
\end{theorem}
\begin{proof}
The proof is along the lines of \cite[Theorem 5.52]{LV}.
\end{proof}

Theorem \ref{ac} will be needed in 
equation~\eqref{ineq11} below.
This is 
one reason why we require (\ref{ineqDCN}) to hold for all $U \in \DC_N$,
as opposed to just $U_N$. (We note that the distinction between
these two definitions disappears in nonbranching spaces, as will be
shown elsewhere.)

\section{Sobolev inequality and global Poincar\'e inequality} \label{global}

\begin{definition}
Given $f \in \Lip(X)$, put
\begin{equation} \label{nablamoins}
|\nabla^- f|(x) \: = \: \limsup_{y\to x} \frac{[f(y)-f(x)]_-}{d(x,y)} \: = \:
\limsup_{y\to x} \frac{[f(x)-f(y)]_+}{d(x,y)}.
\end{equation}
\end{definition}

Here $a_+ \: = \: \max(a, 0)$ and $a_- \: = \: \max(-a, 0)$.
Note that $|\nabla^- f|(x) \: \le \: |\nabla f|(x)$, the latter
being defined as in (\ref{ug}).

\begin{theorem} \label{Sobmess}
Given $N \in (1, \infty)$ and $K > 0$, suppose that $(X,d,\nu)$ has $N$-Ricci
curvature bounded below by $K$.
Then for any positive Lipschitz function $\rho_0 \in \Lip(X)$ with
$\int_X \rho_0 \: d\nu \: = \: 1$, one has
\begin{equation}
N - N \int_X \rho_0^{1 - \frac{1}{N}} \: d\nu \: \le \:
\int_X \theta^{(N,K)}(\rho_0, |\nabla^- \rho_0|) \: d\nu,
\end{equation}
where 
\begin{align} \label{theta}
\theta^{(N,K)}(r,g) \: = \: r \: \sup_{\alpha \in [0, \pi]} \left[
\frac{N-1}{N} \: \frac{g}{r^{1+\frac{1}{N}}} \: \right. 
& \sqrt{\frac{N-1}{K}} \: \alpha \: + \:
N \left( 1 - \left( \frac{\alpha}{\sin(\alpha)} \right)^{1-\frac{1}{N}}
\right) \: + \\
& \left.  (N-1) \: \left( \frac{\alpha}{\tan(\alpha)} - 1\right)
r^{-\frac{1}{N}}
\right]. \notag
\end{align}
\end{theorem}
\begin{proof}
We recall the definitions of $U_N$ and $H_{N,\nu}$ from
(\ref{UN}) and (\ref{HN}).
Applying Definition \ref{N-Riccidef} with $U \: = \: U_N$,
any two probability measures $\mu_0 \: = \: \rho_0 \, \nu$
and $\mu_1 \: = \: \rho_1 \, \nu$ can be joined by a 
Wasserstein geodesic $\{\mu_t\}_{t \in [0,1]}$, arising from an
optimal dynamical transference plan,
along which the following inequality holds :
\begin{align} \label{CDdef}
& H_{N,\nu}(\mu_t) \: \le \\
& N \: - \: N \: 
\int_{X \times X} \left[
\tau^{(1-t)}_{K, N}(d(x_0, x_1)) \cdot \rho_0^{- 
\frac{1}{N}}(x_0)
\: + \: 
\tau^{(t)}_{K, N}(d(x_0, x_1)) \cdot \rho_1^{- \frac{1}{N}}(x_1)
\right] \: d\pi(x_0, x_1).\notag
\end{align}
By Theorem \ref{ac},
$\mu_t$ is absolutely continuous with respect
to $\nu$.

Given a positive function $\rho_0 \in \Lip(X)$, put 
$\mu_0 \: = \: \rho_0 \: \nu$ and
$\mu_1 \: = \: \nu$. Put $\phi(t) \: = \: H_{N,\nu}(\mu_t)$.
In the proof of \cite[Proposition 3.36]{LV} it was shown that
\begin{align} \label{ineq11}
- \: \limsup_{t \rightarrow 0} \frac{\phi(t) \: - \: \phi(0)}{t} 
\: & \le \:
\int_\Gamma U_N^{\prime \prime}(\rho(\gamma(0))) \:
|\nabla^-\rho_0|(\gamma(0)) \: d(\gamma(0), \gamma(1)) \: d\Pi(\gamma)
\\
& = \:
\frac{N-1}{N} \int_X 
\frac{|\nabla^-\rho_0|(x_0)}{\rho_0(x_0)^{1+\frac{1}{N}}} 
\: d(x_0, x_1) \: d\pi(x_0, x_1). \notag
\end{align}
On the other hand,
from (\ref{CDdef}),
\begin{equation}
\phi(t) \: \le \: 
N \: - \: N \: 
\int_{X \times X} \left[
\tau^{(1-t)}_{K, N}(d(x_0, x_1)) \cdot \rho_0^{- 
\frac{1}{N}}(x_0)
\: + \: 
\tau^{(t)}_{K, N}(d(x_0, x_1))
\right] \: d\pi(x_0, x_1)
\end{equation}
and so
\begin{equation} \label{andso}
\frac{\phi(t) \: - \: \phi(0)}{t} \: \le \: 
- \: N \: 
\int_{X \times X} \left[
\frac{\tau^{(1-t)}_{K, N}(d(x_0, x_1))-1}{t} \cdot \rho_0^{- 
\frac{1}{N}}(x_0)
\: + \: 
\frac{\tau^{(t)}_{K, N}(d(x_0, x_1))}{t}
\right] \: d\pi(x_0, x_1).
\end{equation}
Then
\begin{align} \label{ineq2}
& \limsup_{t \rightarrow 0} \frac{\phi(t) \: - \: \phi(0)}{t} \:
\le \: - \: N \: \int_{X \times X}
\left(  \frac{ \sqrt{\frac{K}{N-1}} \: d(x_0,x_1)}{
\sin \left( \sqrt{\frac{K}{N-1}} d(x_0,x_1) \right)} \right)^{1 - \frac{1}{N}} 
\: d\pi(x_0, x_1) \: + \\
& \int_{X \times X}
\left[ 1 \: + \: (N-1) \: \sqrt{\frac{K}{N-1}} \: d(x_0,x_1) \:
\cot \left( \sqrt{\frac{K}{N-1}} \: d(x_0,x_1) \right) \right] \: \rho_0^{-
\frac{1}{N}}(x_0) \: d\pi(x_0, x_1). \notag
\end{align}
Combining (\ref{ineq11}) and (\ref{ineq2}), and slightly rewriting the
result, gives
\begin{align} \label{combine2}
& - \: 
\frac{N-1}{N} \int_X 
\frac{|\nabla^-\rho_0|(x_0)}{\rho_0(x_0)^{1+\frac{1}{N}}} 
\: d(x_0, x_1) \: d\pi(x_0, x_1) \: \le \\ 
& N \: \int_{X \times X} \left[1 \: - \:
\left(  \frac{ \sqrt{\frac{K}{N-1}} \: d(x_0,x_1)}{
\sin \left( \sqrt{\frac{K}{N-1}} d(x_0,x_1) \right)} \right)^{1 - \frac{1}{N}} 
\right]
\: d\pi(x_0, x_1) \: + \notag \\
& (N-1) \: \int_{X \times X}
\left[ \sqrt{\frac{K}{N-1}} \: d(x_0,x_1) \:
\cot \left( \sqrt{\frac{K}{N-1}} \: d(x_0,x_1) \right) \: - \: 1 \right]
\: \rho_0^{-
\frac{1}{N}}(x_0) \: d\pi(x_0, x_1) \: -  \notag \\
& H_{N,\nu}(\mu), \notag
\end{align}
or
\begin{align} 
& H_{N,\nu}(\mu) \: \le \: 
\frac{N-1}{N} \int_X 
\frac{|\nabla^-\rho_0|(x_0)}{\rho_0(x_0)^{1+\frac{1}{N}}} 
\: d(x_0, x_1) \: d\pi(x_0, x_1) \: + \\
& N \: \int_{X \times X} \left[1 \: - \:
\left(  \frac{ \sqrt{\frac{K}{N-1}} \: d(x_0,x_1)}{
\sin \left( \sqrt{\frac{K}{N-1}} d(x_0,x_1) \right)} \right)^{1 - \frac{1}{N}} 
\right]
\: d\pi(x_0, x_1) \: + \notag \\
& (N-1) \: \int_{X \times X}
\left[ \sqrt{\frac{K}{N-1}} \: d(x_0,x_1) \:
\cot \left( \sqrt{\frac{K}{N-1}} \: d(x_0,x_1) \right) \: - \: 1 \right]
\: \rho_0^{-
\frac{1}{N}}(x_0) \: d\pi(x_0, x_1). \notag
\end{align}
Replacing $\sqrt{\frac{K}{N-1}} \: d(x_0,x_1)$ by $\alpha$, we get
only a weaker inequality by taking the sup over $\alpha \in 
[0, \pi]$. The theorem follows. 
\end{proof}

%\begin{remark}
%In the proof of Theorem \ref{Sobmess}, we use Theorem \ref{ac}
%in (\ref{ineq11}).
%\end{remark}

In order to clarify the nature of the inequality of Theorem
\ref{Sobmess}, we derive a slightly weaker inequality.  First, we
prove an elementary estimate.

\begin{lemma} \label{concavelem}
For $x \in [0, \pi]$, one has
\begin{equation} \label{toshow1}
\frac{x}{\tan(x)} \: \le \: 1 - \frac{x^2}{3}
\end{equation}
and
\begin{equation} \label{toshow2}
1 - \left( \frac{x}{\sin(x)} \right)^{1-\frac{1}{N}} \: \le \:  - 
\left( 1 - \frac{1}{N} \right)\frac{x^2}{6}.
\end{equation}
\end{lemma}
\begin{proof}
Put
\begin{equation}
F(x) \: = \: x \: - \: \frac{\sin(x) \cos(x)}{1 - \frac23 \sin^2(x)}.
\end{equation}
Then
\begin{equation}
F^\prime(x) \: = \: \frac49 \: \frac{\sin^4(x)}{\left( 
1 - \frac23 \sin^2(x) \right)^2} \: \ge \: 0.
\end{equation}
As $F(0) = 0$, it follows that $F(x) \ge 0$ for $x \in [0, \pi]$, so
\begin{equation} \label{tech}
x \left( 1 - \frac23 \sin^2(x) \right) \: \ge \: \sin(x) \: \cos(x).
\end{equation}
Putting 
\begin{equation}
G(x)  \: = \: \frac{x}{\tan(x)} \: + \: \frac13 \: x^2
\end{equation}
and using (\ref{tech}), one obtains
\begin{equation}
G^\prime(x) \: = \: - \: \frac{x}{\sin^2(x)} \: + \: \frac{1}{\tan(x)}
\: + \: \frac23 \: x \: \le \: 0.
\end{equation}
As $G(0) =1$, we have
\begin{equation}
\frac{x}{\tan(x)} \: + \: \frac13 \: x^2 \: \le \: 1,
\end{equation}
which proves (\ref{toshow1}).

Next, from (\ref{toshow1}) we have
\begin{equation}
\frac{1}{x} \: - \: \frac{1}{\tan(x)} \: \ge \: \frac{x}{3}. 
\end{equation}
Integrating gives
\begin{equation}
\ln \left( \frac{x}{\sin(x)} \right) \: \ge \: \frac{x^2}{6},
\end{equation}
so
\begin{equation}
\frac{x}{\sin(x)} \: \ge \: e^{\frac{x^2}{6}} \: \ge \:
\left( 1 + \left( 1 - \frac{1}{N} \right) \frac{x^2}{6} 
\right)^{\frac{1}{1-\frac{1}{N}}}.
\end{equation}
Thus
\begin{equation}
\left( \frac{x}{\sin(x)} \right)^{{1-\frac{1}{N}}} \: \ge \: 
1 + \left( 1 - \frac{1}{N} \right) \frac{x^2}{6}.
\end{equation}
This proves (\ref{toshow2}).
\end{proof}

We now prove a Sobolev-type inequality.
\begin{theorem} \label{Sobineq}
Given $N \in (1, \infty)$ and $K > 0$, suppose that $(X,d,\nu)$ has $N$-Ricci
curvature bounded below by $K$.
Then for any nonnegative Lipschitz function $\rho_0 \in \Lip(X)$ with
$\int_X \rho_0 \: d\nu \: = \: 1$, one has
\begin{equation} \label{Sobolev}
N - N \int_X \rho_0^{1-\frac{1}{N}} \: d\nu \: \le \:
\frac{1}{2K} \left( \frac{N-1}{N} \right)^2 \int_X
\frac{\rho_0^{-1-\frac{2}{N}}}{\frac13 + \frac23 \rho_0^{-\frac{1}{N}}} \:
|\nabla^- \rho_0|^2 \: d\nu.
\end{equation}
\end{theorem}
\begin{proof}
If $\rho_0$ is positive then
using Lemma \ref{concavelem}, we can estimate the function 
$\theta^{(N,K)}(r,g)$ of (\ref{theta}) by
\begin{align}
\theta^{(N,K)}(r,g) \: & \le \: 
r \: \sup_{\alpha \in [0, \pi]} 
\left[
\frac{N-1}{N} \: \frac{g}{r^{1+\frac{1}{N}}} \:
\sqrt{\frac{N-1}{K}} \: \alpha \: - \: \frac{N-1}{6} \: \alpha^2 \:
\left( 1 \: + \: 2 \: r^{-\frac{1}{N}} \right)
\right] \\
& \le \: \frac{1}{2K} \left( \frac{N-1}{N} \right)^2 
\frac{r^{-1-\frac{2}{N}}}{\frac13 + \frac23 r^{-\frac{1}{N}}} \:
g^2. \notag
\end{align}
The theorem in this case follows from Theorem \ref{Sobmess}.
The case when $\rho_0$ is nonnegative can be handled by
approximation with positive functions.
\end{proof}

To put Theorem \ref{Sobineq} into a more conventional form,
we prove a slightly weaker inequality.

\begin{theorem} \label{Sobineq2}
Given $N \in (2, \infty)$ and $K > 0$, suppose that $(X,d,\nu)$ has $N$-Ricci
curvature bounded below by $K$.
Then for any nonnegative Lipschitz function $f \in \Lip(X)$ with
$\int_X f^{\frac{2N}{N-2}} \: d\nu \: = \: 1$, one has
\begin{equation} \label{Sobolev2}
1 -  \left( \int_X f \: d\nu \right)^{\frac{2}{N+2}} \: \le \:
\frac{6}{KN} \left( \frac{N}{N-2} \right)^2 \int_X
|\nabla^- f|^2 \: d\nu.
\end{equation}
\end{theorem}
\begin{proof}
Put $\rho_0 \: = \: f^{\frac{2N}{N-2}}$. 
From (\ref{Sobolev}) we have
\begin{equation}
N - N \int_X \rho_0^{1-\frac{1}{N}} \: d\nu \: \le \:
\frac{3}{2K} \left( \frac{N-1}{N} \right)^2 \int_X
\rho_0^{-1-\frac{2}{N}} \:
|\nabla^- \rho_0|^2 \: d\nu,
\end{equation}
which gives
\begin{equation}
1 - \int_X f^{2\left(\frac{N-1}{N-2}\right)} \: d\nu \: \le \:
\frac{6}{KN} \left( \frac{N-1}{N-2} \right)^2 \int_X
|\nabla^- f|^2 \: d\nu.
\end{equation}
By H\"older's inequality,
\begin{equation}
\int_X f^{2\left(\frac{N-1}{N-2}\right)} \: d\nu \: \le \:
\left( \int_X f^{\frac{2N}{N-2}} \: d\nu \right)^{\frac{N}{N+2}} \:
\left( \int_X f \: d\nu \right)^{\frac{2}{N+2}}.
\end{equation}
The theorem follows.
\end{proof}

Putting (\ref{Sobolev2}) into a homogeneous form reveals
the content of Theorem \ref{Sobineq2}: there
is a bound of the form
$\parallel f \parallel_{\frac{2N}{N-2}} \: \le
F \left( \parallel f \parallel_1,
\parallel \nabla^- f \parallel_2 \right)$ for some
appropriate function $F$. This is
an example of Sobolev embedding. Of course
equation (\ref{Sobolev2}) is not sharp, due to the
many approximations made.

Finally, we prove a sharp global Poincar\'e inequality.

\begin{theorem} \label{globalP}
Given $N \in  (1, \infty)$ and $K > 0$, 
suppose that $(X,d, \nu)$ has $N$-Ricci curvature
bounded below by $K$. Suppose that $f \in \Lip(X)$ has
$\int_X f \: d\nu \: = \: 0$. Then
\begin{equation}
\int_X f^2\,d\nu \leq 
\frac{N-1}{KN} \int_X |\nabla^- f|^2\,d\nu. 
\end{equation}
\end{theorem}
\begin{proof}
Without loss of generality we may assume that
$\max |f| \: \le \: 1$.
Given $\epsilon \in (-1, 1)$,
put $\rho_0 \: = \: 1 \: + \: \epsilon \: f$.
Then $\rho_0 > 0$ and $\int_X \rho_0 \: d\nu \: = \: 1$.
For small $\epsilon$,
\begin{equation} \label{end3}
N - N \int_X \rho_0^{1-\frac{1}{N}} \: d\nu \: = \: 
\epsilon^2 \: \frac{N-1}{2N} \int_X f^2 \: d\nu \: + \: 
O(\epsilon^3)
\end{equation}
and
\begin{equation}
\frac{1}{2K} \left( \frac{N-1}{N} \right)^2 \int_X
\frac{\rho_0^{-1-\frac{2}{N}}}{\frac13 + \frac23 \rho_0^{-\frac{1}{N}}} \:
|\nabla^- \rho_0|^2 \: d\nu \: = \:
\frac{\epsilon^2}{2K} \left( \frac{N-1}{N} \right)^2 \int_X
|\nabla^- f|^2 \: d\nu \: + \: 
O(\epsilon^3).
\end{equation}
Then the result follows from Theorem \ref{Sobineq}.
\end{proof}

\begin{remark}
1. In the case of an $N$-dimensional Riemannian manifold with 
$\Ric \: \ge \: K \: g$,
one recovers the Lichnerowicz inequality for the lowest positive
eigenvalue of the Laplacian
\cite{Lichnerowicz}. It is sharp on round spheres. \\
2. The case $N = \infty$ was treated by similar means
in \cite[Theorem 6.18]{LV}.
\end{remark}

\appendix

\section*{Appendix: Stability of dynamical transference plans} 

\def\thesection{A}

In this appendix we prove a general compactness theorem for probability
measures on geodesic paths. This theorem is used
to show that the condition $\Cdem$ is preserved under
measured Gromov-Hausdorff limits.

\begin{lemma} \label{lemma1} Let $X$ be a compact length space. Given 
$\epsilon > 0$, there is a $\delta > 0$ with the following property.  Suppose that
$Y$ is a compact length space and $f \: : \: Y \rightarrow X$ is a $\delta$-approximation.
Let $\gamma \: : \: [0,1] \rightarrow Y$ be a 
geodesic.
Then there is a geodesic $T(\gamma) \: : \: [0,1] \rightarrow X$ so that for all
$t \in [0,1]$,
$d_X(T(\gamma)(t), f(\gamma(t))) \: \le \: \epsilon$.
\end{lemma}
\begin{proof}
Suppose that the lemma is not true.  Then there is some $\epsilon > 0$ along with \\
1. A sequence of compact metric spaces $\{Y_i\}_{i=1}^\infty$, \\
2. $\frac{1}{i}$-approximations $f_i \: : \: Y_i \rightarrow X$ and\\
3.  Geodesics $\gamma_i \: : \: [0,1] \rightarrow Y_i$\\
so that for each geodesic $\gamma^\prime \: : \: [0,1] \rightarrow X$, there is some
$t_{i, \gamma^\prime} \in [0,1]$ with 
\begin{equation} \label{eqn1}
d_X(\gamma^\prime(t_{i, \gamma^\prime}), f_i(\gamma_i(t_{i, \gamma^\prime}))) \: > \: \epsilon.
\end{equation}

After passing to a subsequence, we may assume that 
$\{ f_i \circ \gamma_i\}_{i=1}^\infty$ converges uniformly to a geodesic
$\gamma_\infty \: : \: [0,1] \rightarrow X$. After passing to a further subsequence, we
may assume that $\lim_{i \rightarrow \infty} t_{i, \gamma_\infty} \: = \: t_\infty$ for some
$t_\infty \in [0,1]$. Then
\begin{align}
d_X(\gamma_\infty(t_{i, \gamma_\infty}), f_i(\gamma_i(t_{i, \gamma_\infty})))
\: & \le \:
d_X(\gamma_\infty(t_{i, \gamma_\infty}), \gamma_\infty(t_{\infty})) \: + \:
d_X(\gamma_\infty(t_\infty), f_i(\gamma_i(t_\infty)) \: + \\
& \: \: \: \: \: \: d_X(f_i(\gamma_i(t_\infty), f_i(\gamma_i(t_{i, \gamma_\infty}))) \notag \\
& \le \: \diam(X) \: |t_{i, \gamma_\infty}-t_{\infty}| \: + \: 
d_X(\gamma_\infty(t_\infty), f_i(\gamma_i(t_\infty)) \: + \notag \\
& \: \: \: \: \: \: \frac{1}{i} \: + \: \diam(Y_i) \: |t_{i, \gamma_\infty}-t_{\infty}|. \notag
\end{align}
Then the right-hand side converges to~0 as $i\to\infty$, which 
contradicts (\ref{eqn1}) with $\gamma^\prime \: = \: \gamma_\infty$.
\end{proof}

\begin{lemma}
One can choose the map $T$ in Lemma \ref{lemma1} to be a measurable map
from $\Gamma(Y)$ to $\Gamma(X)$.
\end{lemma}
\begin{proof} 
This follows from \cite[Theorem A.5]{Zimmer}, as
\begin{equation}
\Bigl\{ (\gamma_1, \gamma_2) \in \Gamma(X) \times \Gamma(Y) \: : \:
\text{ for all } t \in [0,1], \: 
d_X(\gamma_1(t), f(\gamma_2(t))) \: \le \: \epsilon\Bigr\}
\end{equation}
is a Borel subset of $\Gamma(X) \times \Gamma(Y)$ and for
each $\gamma_2 \in \Gamma(Y)$,
\begin{equation}
\Bigl\{ \gamma_1 \in \Gamma(X) \: : \:
\text{ for all } t \in [0,1], \: 
d_X(\gamma_1(t), f(\gamma_2(t))) \: \le \: \epsilon\Bigr\}
\end{equation}
is compact.
\end{proof}

\begin{theorem} \label{thm1}
Let $\{(X_i,d_i)\}_{i=1}^\infty$ 
be a sequence of compact length spaces that converges in the
Gromov-Hausdorff topology to a compact length space $(X,d)$.
Let $f_i:X_i\to X$ be $\var_i$-approximations, with
$\var_i\to 0$, that realize the
Gromov-Hausdorff convergence. 
For each $i$, let $\Pi_i$ be a Borel probability measure
on $\Gamma(X_i)$. Let $\pi_i$ and $\{\mu_{i,t}\}_{t \in [0,1]}$ be 
the associated transference plan and measure-valued path. 
Then after passing to a subsequence, there is a dynamical
transference plan $\Pi$ on $X$, with associated transference plan $\pi$,
and measure-valued path $\{\mu_t\}_{t\in [0.1]}$, such that\\
(i) $\lim_{i \rightarrow \infty} (f_i,f_i)_*\pi_i \: = \: \pi$ in the
weak-$*$ topology on $P(X\times X)$;\\
(ii) $\lim_{i \rightarrow \infty} (f_i)_*\mu_{i,t} \: = \: \mu_t$ for all $t \in [0,1]$.
\end{theorem}
\begin{proof}
Let $T_i\: : \: \Gamma(X_i)\rightarrow \Gamma(X)$ 
be the map constructed by means of Lemma \ref{lemma1}
with $\var=\var_i$, $Y=X_i$, $f=f_i$.
After passing to a convergent subsequence, we can assume that
$\lim_{i \rightarrow \infty} (T_i)_* \Pi_i \: = \: \Pi$ in the weak-$*$ topology, for some
$\Pi \in P(\Gamma(X))$. Given $F \in C(X \times X)$, we have
\begin{equation}
\int_{X \times X} F \: d\pi \: = \:
\int_{\Gamma(X)} F(\gamma(0), \gamma(1)) \: d\Pi(\gamma) \: = \:
\lim_{i \rightarrow \infty} \int_{\Gamma(X_i)} F\bigl(T_i(\gamma_i)(0), 
T_i(\gamma_i)(1)\bigr)
\: d\Pi_i(\gamma_i).
\end{equation}
By Lemma \ref{lemma1} and the uniform continuity of $F$,
\begin{align}
\lim_{i \rightarrow \infty} \int_{\Gamma(X_i)} F\bigl(T_i(\gamma_i)(0), T_i(\gamma_i)(1)\bigr)
\: d\Pi_i(\gamma_i) \: & = \:
\lim_{i \rightarrow \infty} \int_{\Gamma(X_i)} F(f_i(\gamma_i(0)), f_i(\gamma_i(1)))
\: d\Pi_i(\gamma_i) \\ 
& = \: \lim_{i \rightarrow \infty} \int_{X_i \times X_i} 
F\bigl(f_i(x_i), f_i(x_i^\prime)\bigr) \: d\pi_i(x_i, x_i^\prime) \notag \\ 
& = \: \lim_{i \rightarrow \infty} \int_{X \times X} F \: d(f_i, f_i)_* \pi_i. \notag
\end{align}
This proves (i). Similarly, for $t \in [0,1]$ and $F \in C(X)$,
\begin{equation}
\int_{X} F \: d\mu_t \: = \:
\int_{\Gamma(X)} F(\gamma(t)) \: d\Pi(\gamma) \: = \:
\lim_{i \rightarrow \infty} \int_{\Gamma(X_i)} F(T_i(\gamma_i)(t))
\: d\Pi_i(\gamma_i).
\end{equation}
By Lemma \ref{lemma1} and the uniform continuity of $F$,
\begin{align}
\lim_{i \rightarrow \infty} \int_{\Gamma(X_i)} F(T_i(\gamma_i)(t))
\: d\Pi_i(\gamma_i) \: & = \:
\lim_{i \rightarrow \infty} \int_{\Gamma(X_i)} F(f_i(\gamma_i(t)))
\: d\Pi_i(\gamma_i) \\ 
& = \: \lim_{i \rightarrow \infty} \int_{X_i} 
F(f_i(x_i)) \: d\mu_{i,t}(x_i) \notag \\ 
& = \: \lim_{i \rightarrow \infty} \int_{X} F \: d(f_i)_* \mu_{i,t}. \notag
\end{align}
This proves (ii).
\end{proof}

For reference we give a slight variation of Lemma \ref{lemma1},
although it is not needed in the body of the paper.

\begin{lemma} Let $X$ be a compact length space. Choose points
$x, x^\prime \in X$. Given 
$\epsilon > 0$, there is a $\delta = \delta(x, x^\prime)  > 0$ with the following property.  Suppose that
$Y$ is a compact length space and $f \: : \: Y \rightarrow X$ is a $\delta$-approximation.
Given $y \in f^{-1}(x)$ and $y^\prime \in f^{-1}(x^\prime)$, let $\gamma \: : \: [0,1] \rightarrow Y$ be a 
geodesic joining them.
Then there is a geodesic $T(\gamma) \: : \: [0,1] \rightarrow X$ from $x$ to $x^\prime$ so that for all
$t \in [0,1]$,
$d_X(T(\gamma)(t), f(\gamma(t))) \: \le \: \epsilon$.
\end{lemma}
\begin{proof}
The proof is along the same lines as that of Lemma \ref{lemma1}.
\end{proof}
\begin{remark}
In general, one cannot take $\delta$ to be independent of $x$ and $x^\prime$.
\end{remark}

\bibliographystyle{acm}

\end{document}